\documentclass{amsart}%
\usepackage{amsfonts}
\usepackage{amsmath}
\usepackage{amssymb}
\usepackage{graphicx}%
\newtheorem{theorem}{Theorem}
\theoremstyle{plain}

\newtheorem{corollary}{Corollary}

\newtheorem{lemma}{Lemma}

\theoremstyle{definition}

\theoremstyle{remark}

\newtheorem{remark}{Remark}

\numberwithin{equation}{section}

\renewenvironment{smallmatrix}{\begin{matrix}}{\end{matrix}}
\newcommand{\R}{{\mathbb R}}
\def\no{\noindent} 

\newtheorem{thm}{Theorem}[section]

\newtheorem{prop}[thm]{Proposition}
\theoremstyle{definition}

\theoremstyle{remark}

\numberwithin{equation}{section}

\begin{document}

\title[Highly Degenerate harmonic mean curvature flow]{ Highly Degenerate\\
harmonic mean curvature flow}
\author[M.C.  Caputo]{M.C.  Caputo$^{*}$}
\address{Department of Mathematics, University of Texas at Austin, TX}
\email{caputo@math.utexas.edu}
\thanks{$*$: Partially supported by the NSF grant  DMS-03-54639 }
\author[P. Daskalopoulos]{P. Daskalopoulos$^{**}$}
\address{Department of Mathematics, Columbia University, NY}
\email{pdaskalo@math.columbia.edu}
\thanks{$**$: Partially supported by the NSF grants DMS-01-02252, DMS-03-54639  and the EPSRC in the UK }
% ----------------------------------------------------------------
\begin{abstract}
We study the evolution of a weakly convex surface $\Sigma_0$ in $\R^3$  with flat sides by the  Harmonic Mean
Curvature flow. We establish the short time
existence  as well as the  optimal regularity of the surface  and we  show that the boundaries of the flat sides evolve by the curve shortening flow. 
It follows from our results that a weakly convex surface with flat sides of class 
$C^{k,\gamma}$, for some $k\in \mathbb{N}$ and  $0 < \gamma \leq 1$,   remains in the same class
under the flow. This distinguishes this flow
from other, previously studied, degenerate parabolic equations, including  the porous medium equation and the Gauss curvature flow with flat sides, where  the regularity of the solution for $t >0$ does not depend on the regularity of the initial data. 
\end{abstract}

\maketitle
% ----------------------------------------------------------------
\section{Introduction}

We consider the motion of a compact, weakly convex two-dimensional
surface $\Sigma_0$ in space $\R^3$ under the 
{\it harmonic mean curvature flow} (HMCF) 
\begin{equation}
\frac {\partial P}{\partial t} = 
\frac{K}{H}\, N\tag{HMCF}
\end{equation}
where  each point $P$ of $\Sigma_0$
moves in the inward normal direction $N$ 
with velocity equal to the  {\it harmonic mean  curvature}  of the surface,
namely the harmonic mean   
$$\frac KH= \frac {\lambda_1\,\lambda_2}{\lambda_1 + \lambda_2}$$
of the two principal curvatures  $\lambda_1$, $\lambda_2$  of the
surface.

%The paper concerns the study of the evolution of a weakly convex
%surface under the {\it harmonic mean curvature flow}.  Differently
%than other flows such as the Gauss curvature flow with a flat side or
%the Porous Medium Equation, a surface evolving by this flow will
%remain of the same class.  Additionally, the motion of the free boundary
%of a surface with a flat side evolving by the {\it Harmonic Mean
%Curvature flow} is specified at each positive time and it moves
%by the curve shortening flow. 

The existence of solutions to the HMCF with strictly convex 
smooth initial data was first shown by Andrews in~\cite{A1}. He also  showed 
that, under the HMCF,  strictly convex,   smooth surfaces converge
to round points in finite time. In~\cite{Di},   Di\"eter established
 the short time existence of solutions to the HMCF  with weakly convex smooth initial data
and mean curvature  $H >0$.  More precisely,  Di\"eter showed that  if at time $t=0$ the surface $\Sigma_0$
satisfies $K\geq 0$ and $H>0$, then there exists a unique   strictly
convex smooth solution $\Sigma_t$  of the HMCF defined on  $0<t< \tau$,  for some
$\tau>0$. By the results of Andrews, this solution exists up
to the time where its enclosed volume becomes zero. However,
the highly degenerate case where the
initial data is weakly convex and both $K$ and $H$  vanish in a region is not studied in~\cite{Di} .

We will consider in this work the evolution of a surface $\Sigma_0$ 
with  flat sides by the HMCF.  The parabolic equation describing the motion
of the surface becomes degenerate at points where  both  curvatures $K$ and
$H$  become zero. Our main objective is to   study the  solvability and optimal regularity of the evolving surface  for $t >0$,  by viewing the flow as a {\it free-boundary} problem.
It will be shown that a surface $\Sigma_0$ of class $C^{k,\gamma}$
with $k\in\mathbb{N}$ and $0<\gamma \leq 1$  at $t=0$,    will remain in the same class
for $t >0$. In addition, we will show that 
the strictly convex parts  of the surface
become instantly $C^\infty$ smooth up to the flat sides and the boundaries  
of the flat sides  evolve by the curve shortening flow.

For simplicity we will assume that the  surface $\Sigma_0$ has only one
flat side, namely  $\Sigma
=\Sigma_1 \cup \Sigma_2$,  with $\Sigma_1$ flat and $\Sigma_2$  strictly convex (both principal curvatures are strictly positive). We may  also assume that $\Sigma_1$ lies on the  $z=0$ 
plane and that $\Sigma_2$ lies above this plane since the equation is invariant under rotation and translation.
Therefore, the lower part of the surface $\Sigma_0$ can be written
as the graph of a function 
$$z=h(x,y)$$
over a compact domain $\varOmega \subset  {\R}^2$ containing
the initial flat side $\Sigma_1$.  Let $\Gamma$ denote the boundary of the flat side $\Sigma_1$. We define
$ g=h^p$, for some $0<p<1 $.
Our main assumption on the initial  surface $\Sigma_0$
is that it satisfies the following {\it non-degeneracy condition~$(\star)$}:

\begin{equation}
 | Dg(P) | \geq \lambda
 \qquad \text{and} \qquad 
g_{\tau\tau}(P)\geq \lambda, \qquad \mbox{for all}\,\,\, P\in\,\Gamma  \tag{$\star$}
\end{equation}
\noindent for some number  $\lambda >0$. Here $\tau$ denotes the
tangential direction to the level sets of $g$ and
$g_{\tau\tau}$ denotes the second order
derivative in this direction.

\noindent Under  the above conditions, our main results show that for $t\in (0,T)$:
\begin{enumerate}
 \item The HMCF admits a solution
$\Sigma_t=(\Sigma_1)_t\cup (\Sigma_2)_t$ of class $C^{k,\gamma}$, for some $k\in \mathbb{N}$ and $0<\gamma \leq 1$ depending on $p$, which
is smooth up to $\Gamma_t =\partial (\Sigma_1)_t $.
\item   $(\Sigma_1)_t$ is flat and its boundary $\Gamma_t$ evolves by the 
curve shortening flow.
\end{enumerate}

\noindent The fact that the solution $\Sigma_t$ remains in the class
$C^{k,\gamma}$ distinguishes this flow from other,
previously studied, degenerate free-boundary problems (such as the
Gauss curvature flow with flat sides, the porous medium equation and
the evolution p-laplacian equation) in which the regularity of the solution 
for $t >0$ does not depend on the regularity of the initial data.

%%$$\frac{\partial f}{\partial
%%t}=\frac{K}{H}\sqrt{1+f_x^2+f_y^2}=$$

%% $$=\frac{f_{xx}\,f_{yy}-f_{xy}^2}{(1+f_x^2+f_y^2)^2}\,\frac{(1+f_x^2+f_y^2)^{3/2}}
%% {(1+f_y^2)\,f_{xx}-2\,f_x\,f_y\,f_{xy}+(1+f_x^2)f_{yy}}\,\sqrt{1+f_x^2+f_y^2}$$
%%\noindent Thus, readily:
%% \begin{equation}\label{eqn:main}\frac{\partial f}{\partial
%% t}=\frac{f_{xx}\,f_{yy}-f_{xy}^2}{\,\,(1+f_y^2)\,f_{xx}-2\,f_x\,f_y\,f_{xy}+(1+f_x^2)\,f_{yy}\,\,}
%% \end{equation}
%% \vskip 0.1 in

We define $\mathfrak S$ to be the class of
weakly convex  compact surfaces $\Sigma_0$ in ${\R}^3$ so that
$\Sigma = \Sigma_1 \cup \Sigma_2$, where $\Sigma_1$ is a surface contained in the plane $z=0$
and $\Sigma_2$ is a strictly convex and smooth surface contained in the half-space $z\geq 0$. The main result states as follows:

%\noindent {\bf Definition.} {\em Given a surface $\Sigma\in \mathfrak S$, we
%say that a family of weakly convex surfaces $\Sigma_t$,
%with $0<t\leq T$, $T>0$, is a solution to the HMCF with initial surface
%$\Sigma$ if:
%\begin{enumerate}
%\item[i.] equation HMCF  is satisfied  at  any point $P\in \Sigma_t$
%such that $H>0$;

%\item[ii.] $\mbox{dist}\, (\Sigma_t,\Sigma\,)\rightarrow\,0, \, \text
%{as}\,\,t\rightarrow 0 $, 
%\end{enumerate}
%where  $\stackrel{\rightarrow}{N}$, 
%$K$ and $H$ are respectively  the outward normal, the Gauss curvature and the Mean curvature of $\Sigma_t$ at the point $P$.}

%
%\section{Main Theorem} \label{sec:shorttime}
\noindent {\bf Main Theorem.}\label{thm1} {\em Assume that at time
$t=0$, $\Sigma_0$ is a weakly convex  compact surface in ${\R}^3$
which belongs to the class $\mathfrak S$ so  that the function $g=h^p$ defined as above  is smooth up to the
interface $\Gamma$ and satisfies ($\star$). Then,
there exists a time $T>0$ such that the HMCF 
admits a solution $\Sigma_t\in \mathfrak S$ on $[0,T)$. Moreover, the function $g=h^p$, defined as above for $\Sigma_t$, is smooth up to the interface $z=0$ and satisfies~($\star$) on $(0,T]$.
In particular, the interface $\Gamma_t$ between the flat side and
the strictly convex side is a smooth curve for all $t$ in $0<t\leq T$ and it evolves by the curve shortening flow.}

%\noindent In the following sketch we underline the main aspects of the proof of the Main Theorem.
\smallskip

\noindent {\it Sketch of the proof}.
A standard  computation shows that
when $\Sigma_t$ solves the HMCF,  the function
$h$ evolves by the equation 
\begin{equation}\label{eqn-hhh}
 h_t =
\frac{h_{zz}h_{yy}-h_{zy}^2 }{(1+h_y^2)h_{zz}-2\,h_z h_y
h_{zy}-(1+h_z^2)h_{yy}}\,\,\,\,\mbox{on}\,\,\,\,\, z>0.
\end{equation}

The HMCF can be seen as a free boundary problem arising from the
degeneracy near the flat side of the fully nonlinear parabolic PDE
which describes the flow. We will  show in section \ref{global} that, 
via a global change of coordinates,   this  free boundary problem is equivalent to an   {\it initial value problem} on $\mathcal D \times [0,T]$, with  $D = \{ (u,v) ; u^2+v^2 \leq 1 \}$,
namely 
\begin{equation}\label{eqn-MMM}
\left\{
\begin{array}{ll}
 Mw =0
& on\,\,{\mathcal{D}}\times [0,T]\\
w=w_0 & at \,\,\,\,\,\,\,t=0%
\end{array}%
\right.
\end{equation} 

\vskip 0.1 in
\noindent The operator $Mw= w_t - F(t,u,v,w,Dw, D^2w)$, is a fully non-linear operator which becomes degenerate
at $\partial \mathcal D$, the boundary of  $\mathcal D$. We will apply the Inverse Function Theorem
between appropriately defined Banach spaces to show that this problem admits
a solution.

The linearization of the operator $M$ at a point $\bar w$ close to the initial data $w_0$,
can be modeled, after straightening the boundary,  on the degenerate equation
\begin{equation}\label{eqn-hl}
f_t =z^2\,  a_{11}  f_{zz} + 2 \, z\, a_{12}  f_{zy} + a_{22}\,  f_{yy} + b_1\, z\,  f_z + b_2 \,  f_y
\end{equation}  
\vskip 0.05 in
\noindent on $z >0$ with no extra conditions
on $f$ along the boundary $z=0$.
We observe that the diffusion in the above equation is
governed by the Riemannian metric $ds^2=d \bar s^2+|dt|$ where
$$d \bar s^2 = {dz^2 \over z^2}+{dy^2}.$$
\vskip 0.01 in
\par\noindent The distance (with respect to the singular metric
$\bar s$) of an interior point ($z >0$) from the boundary ($z=0$) is hence~{\em
infinite}. This fact distinguishes our problem from other, previously
studied, degenerate free-boundary problems such as the degenerate
Gauss curvature flow and the porous medium equation.

\smallskip

{\em The plan of the paper is the following:}  in section \ref{local} we will introduce a local change of coordinates that fixes the
free-boundary $\Gamma$ in equation \eqref{eqn-hhh}.
We will compute the linearization of 
 equation \eqref{eqn-hhh} in this new change of coordinates,  and
 motivate the  use  of the appropriate  Banach spaces $C^{2+\alpha,p}_s$  for our problem.
 The detailed definition of these  Banach spaces will be given in section
 \ref{def:ban}, where we will also present the appropriate 
Schauder estimates  for our problem. In section \ref{deg} we will study the fully-nonlinear
degenerate equations \eqref{eqn-MMM} and establish the short time
existence for such equations in the Banach spaces $C^{2+\alpha,p}_s$.  
The global change of coordinates and the proof of the Main
Theorem will be given in sections \ref{global} and \ref{sec-main} respectively. In  the last section 
we will establish  the comparison principle for  equation \eqref{eqn-hhh}  and characterize our solutions as viscosity solutions.

%\begin{remark} We will actually show that the Main Theorem holds under the weaker condition that $g$ belongs to $C^{2+\alpha}_s$,  
%where the definition of $C^{2+\alpha}_s$ will be given in section 4.
%\end{remark}
{\bf Acknowledgments.} The authors wish to thank G. Huisken for suggesting the problem and R. Hamilton for many stimulating discussions.

\section {Local Change of Coordinates}\label{local}
We will assume throughout this section   that the surface $\Sigma_0$ belongs to the class $\mathfrak S$. Let  $\Sigma_t$ be  a solution to the HMCF on
$[0,T)$, for some $T>0$ in the sense that $\Sigma_t= (\Sigma_1)_t\cup (\Sigma_2)_t$, with $(\Sigma_1)_t$ flat and $(\Sigma_2)_t$ strictly convex. Let $P_0(x_0,y_0,t_0)$ be a point on the interface 
$\Gamma_{t_0}$,  for $t_0>0$ sufficiently small. Then, the strictly
convex part of surface $(\Sigma_{2})_{t_0}$, $t<t_0$ can be expressed locally
around $P_0$ as the graph of a function $z=h(x,y,t)$. 
Let $g$ be
defined by $g=h^p$, for $0<p<1$, such that : g is smooth up to the interface and satisfies condition~($\star$). We can assume,
by rotating the coordinates, that at the point $P_0$ the 
normal vector to $\Gamma_{t_0}$ facing outwards 
the flat side $\Sigma(t_0)$ is
parallel to the $x$-axis,  so that at $P_0$  we have 
$$g_x(P_0) > 0  \qquad \text{and}
\qquad g_y(P_0) =0.$$
Then we solve locally around the point $P_0$ the equation $z=h(x,y,t)$ with respect to $x$. This yields to a map
$x=f(z,y,t)$. The condition on $g$ expressed in terms of $f$ gives the following {\em non-degeneracy condition~($\star \star$)} :
\begin{equation}
 \left (
\begin{split}
-z^{2-p}\,f_{zz}\,\,\, & z^{1-p}\,f_{zy}\\
z^{1-p}\,f_{zy}\,\,\, & \,\, \,\, -f_{yy}
\end{split}\right ) \geq \bar \lambda  \tag {$\star \star$}
\end{equation}
in the sense that both eigenvalues of  the above matrix are bounded from below by a number $\bar \lambda >0$.

%\begin{remark} The Main Theorem can be proven under the weaker assumption that $g\in C^{2+\alpha}_s(\Omega)$, where we
% say that  $g\in C^{2+\alpha}_s(\Omega)$  if and only if $f_l\in C^{2+\alpha,p}_s(\mathcal B_l(P_l))$ for every $l\in I$
%\end{remark}
%We show next some key observations which will play a crucial role in  the proof of the Main Theorem.
 
\noindent Since $f$ is the inverse  of $h$ and the HMCF is invariant
under rotation, the function $f$ satisfies the same equation as $h$ on $z>0$ :
%equation on $z\,>\,0$ :
\begin{equation} \label{eqn-oo}
 f_t =
\frac{f_{zz}f_{yy}-f_{zy}^2 }{(1+f_y^2)f_{zz}-2\,f_z f_y
f_{zy}-(1+f_z^2)f_{yy}}.
\end{equation}
\noindent  We will construct a smooth solution to this equation by using the Inverse Function Theorem. To do so, we will define the Banach space $C^{2+\alpha,p}_s$ in the next
section. According to our notation, the constants $\alpha$ and $p$ indicate ``how the surface becomes
flat'',  while $s$ refers to the hyperbolic metric which governs the
problem.

\noindent We will prove in the next sections  that when  ${\it f} \in C^{2+\alpha,p}_s$ and
satisfies~condition~($\star \star$), then  the equation~(\ref{eqn-oo}) becomes
degenerate at $z=0$ implying that:
\begin{eqnarray} \label{eqn-f} f_t =
\frac{f_{yy}}{1+f_y^2}\,\,\,\mbox{ at the interface}\,\,\, z=0;
\end{eqnarray}
\noindent This is equivalent to say that the free boundary $\Gamma_t$ evolves by
the {\it curve shortening flow}.  

\section{The $C^{2+\alpha,p}_s$ space and Schauder estimates}\label{def:ban}

\noindent Let
$\mathcal{A}$ be a compact subset of the half space $\{\,(z,y)\in \R^2: \,\,
z\ge 0\,\}$ such that $(0,0)\in {\mathcal
A}$. Then, we define:
$$\begin{array}{lcl} \mathcal{A}^{\circ} & := &
\{\,\,y\,\in\,\R\,:(0,y)\,\in\, \mathcal{A}\,\}
\\ \tilde{\mathcal{A}} &: = &
\{(w,y)\,\in \R^2: w= \ln z,\,(z,y)\,\in\,\mathcal{A},\,z \neq\, 0\}\\ Q_T  &: = &
{\mathcal A}\times[0,T],\,\, T>0\\
 Q_T^{\circ}&:= &{\mathcal
A}^{\circ}\times [0,T]\\
\tilde{Q}_T  &: = & \tilde{\mathcal A}\times [0,T].
\end{array}
$$

\noindent Let $0\,<\,p\,<\,1$. Given a function $f$ on
$\mathcal{A}$ we define:
$$\begin{array}{lcl}
f^{\circ}(y) &: = & f(0,y)\\
\tilde{f}(w,y) & : = & e^{-p\,w}\,(\,f(z,y)-f^{\circ}(y)\,)
\end{array}$$

\noindent with \,\,$ w=\ln z,\,\text {for}\,\,z\,>\,0.$

\noindent Analogously, given a function $f$ on $Q_T$ we define:
$$\begin{array}{lcl}
f^{\circ}(y,t)& : = & f(0,y,t)\\
\tilde{f}(w,y,t)& := & e^{-\,p\,w}\,(\,f(z,y,t)-f^{\circ}(y,t)).
\end{array}
$$

%\begin{remark} Note that we choose the simple notation $\tilde{f}$  in spite it depends on $p$. Throughout the paper we will 
%inform the reader of the specific $p$ if necessary.
%\end{remark}

Given a subspace $\mathcal{A}$ as above, we define the 
{\it hyperbolic distance} $\bar s(P_1,P_2)$ between two points 
$P_1=(z_1,y_1)$ and $P_2=(z_2,y_2 ) $ in $\mathcal{A}$  $z_i >0$, $i=1,2$ to be:

\noindent
$$ \bar s(P_1,P_2):= \sqrt{| \ln z_1- \ln z_2|^2+|y_1-y_2|^2},\,\,\,\hbox{  if}\,\,\, 0\,<\, z_1,\,z_2\,\leq 1 $$
\noindent otherwise it is defined to be equivalent to the standard euclidean metric.

We define the 
{\it parabolic hyperbolic distance} between two points $\tilde P_1=(z_1,y_1,t_1)$ and $\tilde P_2=(z_2,y_2,t_2)$  with $ z_i >\,0$,  $i=1,2$ to be:
$$s(\tilde P_1,\tilde P_2):= \bar s(P_1,P_2)+\sqrt{|t_1-t_2|}$$ 
\noindent where $P_1=(z_1,y_1),\,P_2=(z_2,y_2)$.

\vskip 0.1 in
Let $0<\alpha\leq 1$. We define the H\"older space $C^{\alpha,p}_{ \bar s}(\mathcal{A}$) in terms of the above distance.
%.  We say a
%continuous function $f$ on a compact subset $\mathcal{A}$ is H\"older
%continuous with respect to the metric $s$ if there exists $C>0$ such
%that for all points $P_1$ an $P_2$ in $\mathcal{A}$ we have
%$$|\,f(P_1) - f(P_2)\,| \leq C \, s [P_1, P_2]^\alpha$$
\noindent We start defining the H\"older semi-norm:
$$\|\,f\,\|_{H^\alpha_{\bar s}(\mathcal{A})}: = \sup\limits_{P_1\not= P_2\in  \mathcal{A}\,\cap\,\{(x,y)\in \R^2:z\,>\,0\} } \frac{|\,f(P_1) -
f(P_2)\,|}{ s [P_1, P_2]^{\alpha}}\,$$
and the norm
$$\|\,f\,\|_{C^\alpha_{\bar s}(\mathcal{A})}: = \|\,f\,\|_{C^0 (\mathcal{A})}
+ \|\,f\,\|_{H^\alpha_{\bar s}(\mathcal{A})}$$ 
\noindent where $||\,f\,||_{C^{0}({\mathcal{A}})}:=\sup\limits_{P\in\mathcal{A}}|\,f(P)\,|$.

\label{def2} 
\noindent {\it We say that a function $f$ belongs to $C^{\alpha,p}_{\bar  s}(\mathcal{A})$ if $f^{\circ}\in C^{\alpha}(\mathcal{A}^{\circ}) \,\, \text{and}\,\,
\,\tilde{f}\in C^{\alpha}(\tilde{\mathcal{A}})$}. The norm of $f$ in the
space $C^{\alpha,p}_s(\mathcal{A})$ is defined as:
$$||\,f\,||_{C^{\alpha,p}_{\bar s}(\mathcal{A})}:=||\,f^{\circ}\,||_{C^{\alpha}(\mathcal{A}^{\circ})}+||\,\tilde{f}\,||_{C^{\alpha}(\tilde{\mathcal{A}})}.$$
\noindent Moreover, we define: $||\,f\,||_{C^{0,p}(\mathcal
A)}:=||\,f^{\circ}\,||_{C^{0}(\mathcal A^{\circ})}+||\,\tilde
f\,||_{C^{0}(\tilde{\mathcal A})}$\,.
\begin{remark} We observe that $ f(w,y)\in\,C^{\alpha}(\tilde{\mathcal{A}})$ if and only if
 $f(z,y)\,\in\,C^{\alpha}_{\bar s}(\mathcal{A})$, where $w= \ln z$.
\end{remark}

%\begin{remark}
%The hyperbolic metric is weaker than the Euclidean metric, hence, the
%space $C^{\alpha}(\mathcal{A)}$ of the H{\"o}lder functions with
%respect to the Euclidean metric, is a subspace of
%$C^{\alpha}_s(\mathcal{A})$. Also, smooth functions belong to
%$C^{\alpha,p}_s(\mathcal{A})$.
%\end{remark}

 \hspace{1pt} {\it We say that a continuous function $f$ on
$\mathcal{A}$ belongs to $C^{2+p}(\mathcal{A})$ if
$f^{\circ}\,\in\,C^2(\mathcal{A}^{\circ})$ and $f$ has continuous derivatives
 $$ f_z,\,\, f_y,\,\, f_{zz},\,\, f_{zy},\,\, f_{yy}$$
 \noindent in the interior of $\mathcal{A}$, such that

\no
$$z^{-p}\,(f-f^{\circ}),\,\,z^{1-p}\,f_z,\,z^{-p}\,(\,f_y-f_y^{\circ}\,),\,
 z^{2-p}\,f_{zz},\, z^{1-p}\,f_{zy},\,z^{-p}\,(\,f_{yy}-f_{yy}^{\circ}\,)$$
\noindent extend continuously up to the boundary.} The norm of $f$ in the space $C^{2+p}(\mathcal A)$ is defined as follows:
$$||\,f\,||_{C^{2+p}(\mathcal{A})}:=||\displaystyle\sum_{m=0}^2\,D^m_y\,f^{\circ}\,||_{C^{0}(\mathcal{A}^{\circ})}+
\displaystyle\sum_{m+n=0}^2||\,D_z^m\,D_y^n\,\tilde{f}\,||_{C^0(\tilde{\mathcal{A}})}$$

\noindent {\it Given $f\in C^{2+p}(\mathcal{A})$, we say that $f$ belongs
to $C^{2+\alpha,p}_{\bar s}(\mathcal{A})$ if
$$f^{\circ}\in
C^{2+\alpha}(\mathcal{A}^{\circ})\,\,\text{and}\,\, z\, f_z, f_y,\,z^2\,f_{zz},  z\, f_{zy},f_{yy}$$
\vskip 0.1 in
\noindent extend continuously up to the boundary, and the extensions
are H\"older continuous on $\mathcal{A}$ of class
$C_{\bar s}^{{\alpha,p}}(\mathcal{A})$. The norm of $f$ in the space
$C^{2+\alpha,p}_{\bar s}(\mathcal{A})$ is defined as:
$$||\,f\,||_{C^{2+\alpha,p}_{ \bar s}(\mathcal{A})}:=
||\,f^{\circ}\,||_{C^{2+\alpha}(\mathcal{A}^{\circ})}+\displaystyle\sum_{m+n=0}^2||\,z^{m}\,
D_z^m\,D_y^n \,f\,||_{C^{\alpha,p}_{\bar s}(\mathcal{A})}$$}
\begin{remark} It follows by definition that  $\tilde{f}_w=-p\,\tilde{f}+z^{1-p}\,f_z$ and 
$ \tilde{f}_{ww}=-p\,\tilde{f}_z+(1-p)\,z^{2-p}\,f_{zz}$, which implies that:
 $$\displaystyle\sum_{m+n=0}^2||\,z^{m}\,D_z^m\,D_y^n
 \,f\,||_{C^{\alpha,p}_{ \bar s}(\mathcal{A})}\,\simeq \,||\,\tilde{f}\,||_{ C^{2+\alpha}(\tilde{\mathcal{A}})}.$$
\end{remark}

\begin{remark} The function $f\in C^{2+\alpha,p}_{\bar s}(\mathcal{A})\,$ if and only if
$\,f^{\circ}\in C^{2+\alpha}(\mathcal{A}^{\circ})\,$
 and $\,\tilde{f}\in C^{2+\alpha}(\tilde{\mathcal{A}})$, therefore, the following norms are equivalent:
$$||\,f\,||_{C^{2+\alpha,p}_{\bar s}(\mathcal{A})}\simeq \,||\,f^{\circ}\,||_{ C^{2+\alpha}(\mathcal{A}^{\circ})}+\,||\,\tilde{f}\,||_{ C^{2+\alpha}(\tilde{\mathcal{A}})}.$$ 
\end{remark}
\smallskip
\noindent Let $T>0$. The definitions above can be naturally extended
on the space-time domain $Q_T$ by using the parabolic distance
$d{s}^2=d\bar s^2+|dt|$.
We define the space $C^{\alpha}_s({Q_T})$ to be the standard H{\"o}lder space with respect to the metric $d{s}^2$.
{\it We say that a continuous function $f$ on $Q_T$ belongs to $C^{2+p}({Q_T})$ if $f$ has continuous
derivatives
\vskip 0.01 in
 $$ f_t, f_z, f_y, f_{zz}, f_{zy}, f_{yy}$$
\vskip 0.05 in
\noindent in the interior of $Q_T$ and $f^{\circ}$ has continuous
derivatives that extend continuously up to the boundary and
$$z^{-p}\,(f-f^{\circ}),\,\,z^{-p}\,({f}_t-f_t^{\circ}),\,\,z^{1-p}
{f}_z,z^{-p}\,{f}_y,\,\, z^{2-p} {f}_{zz},\,\, z^{1-p}
{f}_{zy},\,\,z^{-p}\,({f}_{yy}-f_{yy}^{\circ})$$
 \noindent extend continuously up to the boundary.}  The norm of $f$ in the space $C^{2+p}({Q_T})$ is defined as follows:

\smallskip
\centerline{$||\,f\,||_{C^{2+p}}:=||\,f^{\circ}\,||_{C^{2}}+\displaystyle\sum_{l+m+2j=0}^2||\, D_z^l\,D_y^m\,D_t^j\,\tilde f ||_{C^{\circ}}$}

\noindent

\noindent {\it The function $f$ belongs to $ C^{2+\alpha,p}_{\bar s}(Q_T)$ if $f\in C^{2+p}(Q_T)$, 
$$ f,f_t,z f_z, f_y\qquad \text{and}\qquad z^2 f_{zz}, z f_{zy}, f_{yy}\,\,\text {belong to}\,\,\, C_{\bar s}^{{\alpha,p}}({Q_T}).$$}

\smallskip
Throughout the paper $k$ will denote a positive integer. We can extend these definitions to
spaces of higher order derivatives. We denote by $C^{k,p}(Q_T)$ the space
of all functions $f$ whose $k$-th order derivatives $
D^i_z\,D^j_y\,D^l_t\,f $, $ i+j+2l= k$ in the interior of $Q_T$ and
$z^{i}\, D^i_z\,D^j_yD^l_t\,( f-f^{\circ}\,)$, $ i+j+2l= k$ exist and belong to the
space $C^{0}(Q_T)$. We define $C^{\infty,p}(Q_T)=\cap_{k}\, C^{k,p}(Q_T)$.

We denote by $C^{k+\alpha,p}_s(Q_T)$ 
the space of all functions $f\in C^{k,p}(Q_T)$ such that
 $z^{i}\, D^i_z\,D^j_y\,D^l_t\, f$, for $ i+j+2l= k$  belong to the space $C^{\alpha,p}_s(Q_T)$.
The space $C^{k+\alpha,p}_s(Q_T)$ is equipped with the norm:
$$||\,f\,||_{C^{k+\alpha,p}_{\bar s}(Q_T)}
: = \sum_{i+j+2l \leq  k}
 ||\,z^{i}\,D^i_zD^j_yD^l_t \,f\,||_{C^{\alpha,p}_{\bar s}(Q_T)}.$$

\begin{remark}\label{rmkfin}
\noindent A function $f\,\in\,C^{k+\alpha,p}_s(Q_T)$ iff
$f^{\circ}\,\in\,C^{k+\alpha}(Q^{\circ}_T)$ and  $\tilde{f}\,\in\,C^{k+\alpha}(\tilde{Q}_T)$. Moreover, 
$$||\,f\,||_{C^{k+\alpha,p}_s(Q_T)} \simeq
||\,f^{\circ}\,||_{C^{k+\alpha,[k/2]+\alpha/2}
(Q^{\circ}_T)}+||\,\tilde{f}\,||_{C^{k+\alpha,[k/2]+\alpha/2}(\tilde{Q}_T)}.$$
\end{remark}
%\begin{Lemma} 

\bigskip
\no In the next paragraph we denote by ${\mathcal S}_0$ the half space $x
\geq 0$ in ${\R}^2$, by $\mathcal S$ the space $\mathcal S =
{\mathcal S}_0 \times [0,\infty)$, and by $\mathcal S_T$ the space
$\mathcal S \times [0,T]$, for $T>0$.
\no The operator $L_k: \, C^{k+2+\alpha,p}_s(Q_T) \to C^{k+\alpha,p}_s(Q_T)$ is defined as:
\begin{equation}\label{eqn:ope}
L_k\,[\,f\,]:= f_t -(\,z^2 a_{11}  f_{zz} + 2 \, z\, a_{12}
 f_{zy} + a_{22}  f_{yy} + b_1 z\,  f_z + b_2 
f_y+c\,f\,)
\end{equation}
\no where the coefficients $\{a_{ij}\}_{i,j}$ are uniformly elliptic and 
$\{a_{ij},\,b_i,\,c\}\,\subseteq\,C^{k+2+\alpha}_s(Q_T)$,
$\{a_{22},\,b_2,\,c\}\subseteq C^{k+2+\alpha,p}_s(Q_T), \,\,i,j=1,2$.

\begin{theorem} (Existence and Uniqueness) \label{thm:exi}
 Let $L_k$ be defined as above. Assume that $\phi \in 
C^{k+\alpha,p}_s(\mathcal S)$ and $f_0 \in 
C^{k+2+\alpha,p}_s(\mathcal S_0)$, and that $\phi$, $f_0$ are compactly
supported in $\mathcal S$ and $\mathcal S_0$, respectively. Then, for
any $T>0$, the initial value problem\,:
\begin{equation}
\label{eqn:cauchy}
\left\{
\begin{array}{cccc}
 L_k\,[ f] &=&\phi
& in\,\,S_T \\
f(\cdot,0)&= &f_0 & on\,\,S_0 
\end{array}
\right.
\end{equation}
admits a unique solution $f \in C^{k+2+\alpha,p}_s(\mathcal
 S_T)$. Moreover
\begin{eqnarray}\label{ineqg} 
||f||_{C^{k+2+\alpha,p}_s(\mathcal S_T)}
 \leq C(T) \left ( ||f_0||_{\mathcal
C^{k+2+\alpha,p}_s(\mathcal S_0)} + ||\phi||_{\mathcal
C^{k+\alpha,p}_s(\mathcal S)} \right )
 \end{eqnarray}
 \no for some constant $C(T)$, depending on $\alpha$, $p$, $k$ and $T$.
\end{theorem}
\vspace{ -0.2 in}
\begin{proof} 
To solve the above Cauchy problem is
equivalent to solve the following Cauchy problems~(\ref{eqn:cauchya}) and (\ref{eqn:cauchyb}).

%It is easy to observe that the interesting case holds when the Lebesgue measure of the supports of $f_0^{\circ}$ and $\phi^{\circ}$ satisfies
%$|(Supp\,f_0)^{\circ}|\geq \eta$, $|(Supp\, \phi)^{\circ}|\geq \eta$,
%for some $\eta\,>0$. 
\noindent The problem~(\ref{eqn:cauchya}) is obtained by evaluating~(\ref{eqn:cauchy}) at $z=0$ and the problem~(\ref{eqn:cauchyb})  is obtained by solving the corresponding one for $\tilde{f}$. 

\begin{equation}\label{eqn:cauchya}
\left\{
\begin{array}{ll}
(L_k)_0\,[ f^{\circ}] =\phi^{\circ}
& in\,\,\R\times[0,T]\\
f^{\circ}(\cdot,0)=(f_0)^{\circ}  & on\,\,\R%
\end{array}%
\right.
\end{equation}
  
\begin{equation}\label{eqn:cauchyb}\left\{
\begin{array}{ll}
 \tilde{L_k}\, [\tilde{f}] =\tilde\phi
& in\,\,S_T\\
\tilde{f}(\cdot,0)=\tilde{f}_0 & on\,\,S_0%
\end{array}%
\right.
\end{equation}
where the operators  $({L_k})_0$ and $\tilde{L_k}$ are defined respectively as follows: 
$$({L_k})_0( f^{\circ})= a_{22}^{\circ}\,f_{yy}^{\circ}+b_2^{\circ}\, f_y^{\circ}+c^{\circ} f^{\circ};$$ 

\begin{equation}\label{eqn:lt}
\tilde L_k\,[\tilde{f}] =\tilde{f}_t -( \hat{a}_{11} \tilde{f}_{ww} +
2\, \hat a_{12} \tilde{f}_{wy} + \hat{ a}_{22} \tilde{f}_{yy} +
\hat b_1 \,\tilde{f}_w+ \hat b_2 \,\tilde{f}_y +\hat c\,\tilde{f}+\hat G)
\end{equation}
\par\noindent
with
$$
\begin{array}{lcl}
\hat a_{ij}(w,y,t) & := &  a_{ij}(x,y,t)\\
\hat b_1(w,y,t) &:= & (2p-1)\,a_{11}(w,y,t)+b_1(x,y,t)\\
 \hat b_2(w,y,t) & := & b_2(x,y,t)
\end{array}
$$
$$
\begin{array}{lcl}
\hat c(w,y,t)& := & e^{-p\,z}[\, p^2\,\hat a_{11}(x,y,t)-2\,p\,\hat
a_{12}(w,y,t)+p\, b_1(x,y,t)\, ]\\

\hat{G}(w,y,t) & := &  \tilde{b_2}(w,y,t)\,g_y^{\circ}(y,t)+\hat a_{22}(w,y,t)
 g_{yy}^{\circ}(y,t). 
\end{array}
$$

By the assumptions on the operator $L_k$ it is clear that the
coefficients of the two operators $({L_k})_0$ and $\tilde{L_k}$ satisfy
classical conditions. We first find the solution $f^{\circ}$
to~(\ref{eqn:cauchya}), then we solve~(\ref{eqn:cauchyb}).  By
classical theory both problems have a unique solution.  \noindent The function $f$ defined by $f(w,y,t):=f^{\circ}(y,t)+z^p\,\tilde{f}(w,y,t)$ is a solution to~(\ref{eqn:cauchy}).  

\phantom{ufveycutyv}

Let $C^{k+\alpha}$ and $C^{k+2+\alpha}$ denote classical  parabolic H\"older spaces. Then the
following inequalities hold:
$$
\begin{small}\begin{array}{rcl}
||f^{\circ}||_{C^{k+2+\alpha}(\R^+\times [0,T])} &\leq & C(T) \left (
||f^{\circ}_0||_{C^{k+2+\alpha}(\R^+)} + ||g^{\circ}||_{\mathcal
C^{k+\alpha}(\R^+)} \right )\\
\\
||\tilde{f}||_{C^{k+2+\alpha}(\tilde{\mathcal S}_T)} &\leq & C(T)
\left ( ||\tilde{f}_0
||_{C^{k+2+\alpha}(\tilde{\mathcal S}_0)} +
||\tilde{\phi}||_{C^{k+\alpha}(\tilde{\mathcal S})} \right)
\end{array}\end{small}
$$
\no It follows that the
solution to~(\ref{eqn:cauchy}) is unique and satisfies the inequality (\ref{ineqg}).
\end{proof}

%Indeed, at $t=0$,
%$f(\cdot,0)=f^{\circ}(\cdot,0)+z^p\,\tilde{f}(\cdot,0)=f_0$. 

%\no {\it Claim:} $L\,f=\phi$. 
%\noindent Indeed:
%$$L\,f=(L\,f)^{\circ}+z^p\,(\widetilde{L\,f})=
%L_0(\,f^{\circ})+z^p\,\tilde{L}(\,\tilde{f})=\phi^{\circ}+z^p\,\tilde{\phi}=\phi$$
%As a consequence, the function $f$ is a solution to~(\ref{eqn:cauchy}). 

\no Next, we define the boxes in which we prove the Schauder estimates. Let $0<r\leq 1$. We denote  by ${\mathcal B}_r(P)$ the box
$${\mathcal B}_r(P) = \left\{\left(\begin{smallmatrix}
z\\y\\t\end{smallmatrix}\right) : \begin{smallmatrix} z\geq 0, |x-x_0|\le e^r\\
|y-y_0|\le r\\ t_0 - r^2\le t \le t_0\end{smallmatrix}\right\}\ $$ around the
point $P=\left(\begin{smallmatrix}
z_0\\y_0\\t_0\end{smallmatrix}\right)$. We set ${\mathcal B}_r$ to be
the box around the point $P=\left(\begin{smallmatrix}
0\\0\\1\end{smallmatrix}\right)$.  
\begin{remark} The choice of the box ${\mathcal B_r}$ is made so that it has the right rescaling. 
The operators  ${L_k}_0$ and $\tilde L_k$  are well understood on the corresponding boxes 
$\mathcal {B}^{\circ}_r$ and ${\tilde{\mathcal B}}_r$.

\end{remark}
\begin{theorem}
\label{thm:sc} (\mbox{ Schauder\, Estimate})
 Assume that all the coefficients
of the operator 
$$\,L\,f = f_t - (z^2\, a_{11} f_{zz} + 2\, z \, a_{12}  f_{zy}
+ a_{22}  f_{yy} +z\, b_1  f_x + b_2  f_y+\,c f)$$
 belong to the space $C^{k+\alpha}_s(\mathcal B_1)$ and that the coefficients $a_{22}, b_2$ and $c$ belong to $ C^{k+\alpha,p}_s(Q_T)$
for some numbers $\alpha$, $p$ in $0<p <1$, $0<\alpha\leq 1$ and satisfy
$$a_{ij} \xi^i \xi^j  \geq \lambda |\xi|^2, \,\,\,
\forall \xi \in {\mathcal R}^2 \setminus \{0 \} \,\,\,\lambda >0$$
$$with \,\,||a_{ij}||_{C^{k+\alpha}_s(Q_T)}, \,||b_i||_{C^{k+\alpha}_s(Q_T)}, \,||a_{22}||_{C^{k+\alpha,p}_s(Q_T)},||b_2||_{C^{k+\alpha,p}_s(Q_T)},
||c||_{C^{k+\alpha,p}_s(Q_T)}
 \leq{ 1\over \lambda}.$$

\noindent Then,
 there exists a constant $C$ depending only on $\alpha$, $\lambda$ and $p$
such that 
$$\|f\|_{ C^{2+ \alpha,p}_s(\mathcal B_{1/2})} \le
C\left(\|f\|_{ C^{0,p} (\mathcal B_1)} + \|L[ f]\|_{ C^{
\alpha,p}_s(\mathcal B_1)}\right)$$ for all functions $f \in
C^{2+\alpha,p}_s(\mathcal B_1)$.
\end{theorem}
\begin{proof}
The proof follows by the same argument as in Theorem~\ref{thm:exi} and classical Schauder  estimates for strictly parabolic operators. 
\end{proof}

\section{The Degenerate Equation on the disc}\label{deg}

\medskip
We will extend in this section the  existence and
uniqueness the Theorem~\ref{thm:exi} to the following class of
linear degenerate equations: 
$$Lw: = w_t  - \, (\, a^{ij} w_{ij} + b^i\, w_{i} + c\,w \, ) $$
on the cylinder ${\mathcal D}  \times [0,T)$, $T >0$,  where
${\mathcal D}$ denotes the unit disk in  ${\R}^2$.
The sub-indices $i,j\in \{x,y \}$ denote differentiation with respect
to the space variables $x,y$.
% and the summation convention is used.
The matrix $\{a^{ij}\}$ is assumed to be symmetric.
Certain assumptions on the coefficients will be made 
so that this class of equations includes, 
under appropriate change of coordinates,
the operator $L$ given by~({\ref{eqn:ope}).

%The function $ \vartheta (x)$  is assumed to be smooth on $\mathcal D$,
%strictly positive in its interior, with 
%$$\vartheta (P) = \text{dist}\, (P, \partial \mathcal D)$$
%for all $P \in \mathcal D$ sufficiently close to $\partial \mathcal D$. 

We define the distance 
function $s$ in ${\mathcal D}$ as follows:
in the interior 
of ${\mathcal D}$, $\bar s$ it is equivalent to
the standard euclidean distance,
while around any boundary point $P \in \partial 
\mathcal D$, $\bar s$ is defined as the pull
back of the distance function induced by the metric
$$d\bar s^2 = \frac {dz^2}{z^2} + dy^2$$
on the half space $\mathcal S_0  = \{ (z,y) :  z \geq 0  \}$,
via a map 
$\varphi: \mathcal S_0  \cap {\mathcal D} \to {\mathcal D}$ that flattens the boundary 
of the disk ${\mathcal D}$ near $P$.

The parabolic distance is defined by
$$s \left [\left(\begin{smallmatrix} P_1\\ t_1\end{smallmatrix}\right),
\left(\begin{smallmatrix} P_2\\ t_2\end{smallmatrix}\right) \right ] =
\bar s(P_1,P_2) + \sqrt {|t_1-t_2|}, \,\,\, P_1,P_2 \in {\mathcal D},\,\,0<t_1\leq t_2$$
We define the spaces $C^{k+\alpha,p}_s({\mathcal D})$ and $C^{k+2+\alpha,p}_s({\mathcal D})$.
\no For a fixed small number $\delta$ in $0 < \delta <1$,
we write  
$${\mathcal D} = {\mathcal D}_{1-\delta/2} \, \cup  ({\displaystyle\bigcup_l}\, 
\left ( {\mathcal D}_{\delta}(P_l) \cap {\mathcal D}) \right )$$
for finite many points $P_l \in \partial {\mathcal D}$, $l \in I$, 
with ${\mathcal D}_{1-\delta/2}$ denoting the disk centered at the origin of radius
$1-\delta/2$ and ${\mathcal D}_{\delta}(P_l)$ denoting the disk of radius
$\delta$ centered at $P_l$. 

We denote by ${\mathcal D}_+$ the half disk
$${\mathcal D}_+ = \{ \, (x,y) \in {\mathcal D} : \,\, x \geq 0 \, \}.$$
We can choose charts $\Upsilon_l : {\mathcal D}_+ \to {\mathcal D}_{\delta}(x_l) \cap 
{\mathcal D}$
which flatten the boundary of ${\mathcal D}$ and such that $\Upsilon_l (0) = P_l$, $ l\in\,I$. Let $\{\psi$, $\psi_l\}$ 
be a partition of unity subordinated to the cover
$$\{ \, {\mathcal D}_{1-\delta/2},\,
( {\mathcal D}_{\delta}(P_l) \cap {\mathcal D}) \,\} $$ of ${\mathcal D}$, with $l\in\, I$.

We define  $C^{k+\alpha,p}_s({\mathcal D})$ 
 to be the space of all functions $w$
on ${\mathcal D}$ such that $w \in C^{k+\alpha} ({\mathcal D}_{1-\delta/2})$ and 
$w \circ \Upsilon_l \in C^{k+\alpha,p}_s({\mathcal D}_+)$
 for all $l \in I$.

Also,  we  define $C^{k+2+\alpha,p}_s({\mathcal D})$ to be the space  of all functions $w$
on ${\mathcal D}$ such that $w \in C^{k+2+\alpha }({\mathcal D}_{1-\delta/2})$ and
$w \circ \psi_l \in C^{k+2+\alpha,p}_s({\mathcal D}_+)$ for all $l \in I$.
Here $C^{k+\alpha}$ and $C^{k+2+\alpha}$ denote the regular H\"older Spaces, 
while $C^{k+\alpha,p}_s ({\mathcal D}_+)$ and $ C^{k+2+\alpha,p}_s({\mathcal D}_+)$ denote the
H\"older Spaces defined in section~\ref{def:ban}.
One can show that both spaces $C^{k+\alpha,p}_s({\mathcal D})$ and $C^{k+2+\alpha,p}_s({\mathcal D})$ are Banach Spaces 
under the  norms
$$ ||w||_{C^{k+\alpha,p}_s({\mathcal D})} = ||\psi \, w||_{C^{k+\alpha} ({\mathcal D}_{1-\delta/2})}
+ \sum_{l} \, ||\psi_l \,( w \circ \Upsilon_l)||_{C^{k+\alpha,p}_s({\mathcal D}_+)}$$ 
and
$$ ||w||_{C^{k+2+\alpha,p}_s({\mathcal D})} = ||\psi \, w||_{C^{k+2+\alpha}(
{\mathcal D}_{1-\delta/2})}
+ \sum_{l} \, ||\psi_l \,( w \circ \Upsilon_l)||_{C^{k+2+\alpha,p}_s({\mathcal D})}.$$
\noindent

The above definitions can be extended in
a straight forward manner to the parabolic spaces 
$C^{\alpha,p}_s(Q)$ and $C^{2+\alpha,p}_s(Q)$ where $Q$ is the 
cylinder $Q= {\mathcal D} \times [0,T]$, for some $T >0$.
\smallskip
\noindent Before we state the main result in this section, we will give the assumptions on the coefficients
of the equation
$$w_t =\, a^{ij} w_{ij} + b^i\, w_{i} + c\,w $$
on the cylinder $Q = {\mathcal D} \times [0,T)$, $i,j=1,2$.

We first assume that for any  $\delta$ in $0<\delta <1$, the  
coefficients $\{a^{ij}\}$, $\{b^i\}$ and $c$
belong to the 
H\"older class $C^{\alpha}(\mathcal D_{1-\delta/2}
 \times [0,T])$, which means that the coefficients 
are of the class $C^{\alpha}$ in the interior of $\mathcal D$.
For a number $\delta$ in $0<\delta <1$,
let  $\Upsilon_l: \mathcal D_+ \to \mathcal D_\delta(P_l) \cap \mathcal D$
be the collection 
of charts which flatten the boundary of $\mathcal D$,
considered above.
We assume that there exists a number $\delta$
so that for every $l \in I$, the coordinate change 
introduced by each of the $\Upsilon_l$ transforms the  
operator 
\begin{equation}\label{opel}
L[w] = w_t - (\,\,  a^{ij}\,  w_{ij} + b^i\, w_{i} + c\,w  \, )\end{equation}
on $\mathcal D_\delta(P_l) \cap \mathcal D$, 
into an operator $\widetilde L_l$ on $\mathcal D_+$ of the form
$$\widetilde L_l \,  [\tilde w] = 
\tilde w_t - (\,  x^2\, \tilde 
a_{11} \, \tilde w_{xx} + 2 \, x \, \tilde a_{12} \, \tilde w_{xy}
+ \tilde a_{22} \, \tilde w_{yy} + x\,\tilde b_1 \,\tilde w_x  + 
\tilde b_2 \, \tilde w_y \,+\tilde c\,\tilde{w})$$
with the coefficients 
$\tilde  a_{ij}$, $\tilde b_i$ and $\tilde c$ 
belonging to the class $C^{k+\alpha}_s(\mathcal D_+)$, with
$a_{22},\,b_2$ and $c\,\in\,C^{k+\alpha,p}_s$ such that:
$$\tilde a_{ij} \xi^i \xi^j  \geq \lambda \, |\xi|^2, 
\qquad \forall \xi \in {\mathcal R}^2 \setminus \{0\}$$ 

\no for some number $\lambda >0$.\\

\no  We need the next Lemma to prove the invertibility of the operator $L$:
\begin{lemma}\label{holder}{(H\"{o}lder Interpolation)}. For every $\epsilon>0$ there exists a constant $C(\epsilon)$ depending on $\epsilon$, $p$, $k$ and $\alpha$ such that for any  $g\in {C^{k+2+\alpha,p}_s (Q_\delta)} $, the following inequality holds:
\begin{equation}\label{eqn:hol}
||\vartheta\,  D\, g\, ||_{C_s^{k+\alpha,p} (Q_\delta)} \leq
\epsilon \, ||\, g\,||_{C_s^{k+2+\alpha,p} (Q_\delta)} + C(\epsilon)
\, ||\,g\,||_{C^{k,p}(Q_\delta)}.
\end{equation}
\noindent where $\vartheta $ behaves like distance to the boundary. 
\end{lemma}
\begin{proof}
It follows by standard arguments. 
\end{proof}

The following existence result readily follows from Theorem~\ref{thm:exi} and the above discussion.

\begin{theorem}\label{thm:exi2} Assume that the operator $L$ satisfies all the above conditions
on the cylinder $Q={\mathcal{D}} \times [0,T]$.
Then, given any function $w_0 \in C^{k+2+ \alpha,p}_s({\mathcal{D}})$ and any function
$g \in C^{k+\alpha,p}_s(Q)$
there exists a unique solution $w \in C^{k+ 2+ \alpha,p}_s (Q_T)$
of the initial value problem
\[
\left\{
\begin{array}{ll}
 Lw =g
& in\,\,Q\\
w(\cdot,0)=w_0 & on\,\,\mathcal D%
\end{array}%
\right.
\]
satisfying
\begin{equation}\label{ineq:www}
||w||_{C^{k+2+\alpha,p}_s(Q)} \leq C(T) \left ( ||w_0||_{C^{k+2+\alpha,p}_s(\mathcal{D})} + ||g||_{C^{k+\alpha,p}_s(Q)} \right )
\end{equation}
The constant $C(T)$ depends  only on the 
numbers $\alpha$, $k$, $\lambda$ and $T$.
\end{theorem}

\begin{proof} 
We can assume, without loss of generality,
that $w_0\equiv 0$ and that 
$g$ is a function in $C_s^{k+\alpha,p}(Q_T)$, which vanishes at $t=0$.
%due to Lemma~\ref{lem:hh}.

For  $\delta >0$, set 
$Q_\delta = \mathcal{D} \times [0,\delta]$ and denote by 
$C_{s,0}^{k+2+\alpha,p}(Q_\delta)$ and
$C_{s,0}^{k+\alpha,p}(Q_\delta)$ the subspaces of
$C_{s}^{k+2+\alpha}(Q_\delta)$ and $C_{s}^{k+\alpha}
(Q_\delta)$ respectively,
consisting out of all functions which vanish identically at $t=0$.
Also, we denote  by $I$ the identity operator on $C_{s,0}^{k+\alpha,p}(Q_\delta)$.
We will show that,  if $\delta$ is sufficiently small,
there exists an  operator $M\, : \, C_{s,0}^{k+\alpha,p}(Q_\delta)
\to C_{s,0}^{k+2+\alpha,p}(Q_\delta)$ such that 
$$ || \, L\,M - I \, || \leq  \frac 12.$$
This implies that the operator 
$L\,M : \, C_{s,0}^{k+\alpha,p}(Q_\delta)
\to C_{s,0}^{k+\alpha,p}(Q_\delta)$ is invertible and therefore 
$L: \,  C_{s,0}^{k+2+\alpha,p}(Q_\delta)
\to C_{s,0}^{k+\alpha,p}(Q_\delta)$ is onto, as desired.
 
We begin by expressing the compact domain $\mathcal D$ as the finite union 

$$\mathcal D = \mathcal D_0 \,\cup{\displaystyle\bigcup_{l\geq 1 }}\, \mathcal D_l$$ 
of compact domains in such a way that  
$$\text{dist}\,(\mathcal D_0, \partial \mathcal D) \geq \frac {\rho}2 >0 $$
and for all $l \geq 1$
$$\mathcal D_l = B_\rho (x_l) \cap \mathcal D$$
with $B_\rho (x_l)$ denoting the ball  
centered at $x_l \in \partial \mathcal D$  of radius 
$\rho >0$. The number $\rho >0$ will 
be determined later.

 The operator $L$ is non-degenerate when restricted on the interior
domain $\mathcal D_0$. Therefore, the classical Schauder theory
for linear parabolic equations implies that $L$ is invertible when restricted on functions which vanish outside $\mathcal D_0$.  
%Notice
%that our H\"older spaces with respect to the hyperbolic metric $s$ on
%the interior domain $\mathcal D _0$ coincide with the standard
%H\"older spaces, where the classical Schauder theory holds true.

We denote by $M_0: \, C_{s,0}^{k+\alpha,p}(\mathcal D_0 \times
[0,\delta]) \to C_{s,0}^{k+2+\alpha,p}(\mathcal D_0 \times
[0,\delta])$ the inverse of the operator $L$ restricted on $\mathcal
D_0$. Next, we consider the domains $\mathcal D_l, \, 
l \geq 1$, close to the boundary of $\mathcal D$,
which can be chosen in such a way that 
the sets $B_{\rho/4}(x_l) \cap \mathcal D $ are disjoint.
Denoting  by $\overline B$ the half unit  ball
$$\overline B= \{(x,y) \in B_1(0) \,; \,\,\, x \geq  0 \, \}$$
and by $\overline Q_\delta$ the cylinder 
$$\overline Q_\delta = \overline B \times [0,\delta]$$
we select  smooth charts
$\Upsilon_l: \overline B \to \mathcal D_l $, which flatten the boundary 
of $\mathcal D$, i.e., they map $\overline B \cap \{ x=0 \} $ 
onto $ \mathcal D_l \cap \partial \mathcal D$ and have
$\Upsilon_l (0) = x_l$. 
This is possible if the number $\rho$ is chosen sufficiently
small. Under the change of coordinates induced by the charts $\Upsilon_l$, the
operator $L$, restricted on each 
$\mathcal D_l \times [0,\delta]$,
 is transformed to  an operator $\bar L_l$ of the form
$$\bar L_l[\bar w ] = \bar w_t - 
(\,x^2 \,\bar a_l^{11}\bar  w_{11} + 2\,x \,\bar a_l^{12}\bar  w_{12}+\bar a_l^{22}\bar  w_{22}+
\,x\bar b_l^1\, \bar w_{1} + \bar b_l^2\, \bar w_{2}+\bar c_l\, \bar w )$$
defined on $\overline B \times [0,\delta]$.
Moreover, the charts $\Upsilon_l$
can be chosen appropriately
so that the coefficients of $\bar L_l$ 
satisfy
$$ \bar a_l^{ij}  \xi_i \xi_j  \geq\bar \lambda |\xi|^2 >0 \qquad \forall \xi
\in \R^2 \setminus \{ 0 \}$$
and
$$||\bar a_l^{ij}||_{C^{k+\alpha,p}_s(\bar Q_\delta)} \quad
||\bar b_l^i||_{C^{k+\alpha,p}_s(\bar Q_\delta)} \quad
||\bar c_l||_{C^{k+\alpha,p}_s(\bar Q_\delta)} \leq  1/\,{\bar 
\lambda}$$

\no for some positive constant $\bar \lambda$.

Each of the operators $ \bar L_l$ has the form of the model operators
previously studied. Denote by ${S}_0$ the half space $x\geq 0$ in
$\R^2$ and by ${S}_\delta$ the space $\mathcal{S}_0 \times
[0,\delta]$.  Also, consider the subspace $\overline
C_{s,0}^{k+\alpha,p}(\mathcal{S}_\delta)$ of
$C_{s,0}^{k+\alpha,p}(\mathcal{S}_\delta)$, consisting out of
functions which are compactly supported on $\mathcal{S}_\delta$.
Then, Theorem~\ref{thm:exi}  implies that for every $l=1,2,...$ there is an
operator $ \overline{M}_l: \overline{C}_{s,0}^{k+\alpha,p}({S}_\delta) \to
{C}_{s,0}^{k+2+\alpha,p}({S}_\delta)$ such that
$$ \bar L_l \,  \overline{ M}_l = I$$
with $I$ denoting the identity operator on 
$\overline C_{s,0}^{k+\alpha,p}(\mathcal{S}_\delta)$.
Let 
$M_l$ be the pull back of the operator $ \overline{ M}_l$  via the chart 
$\Upsilon_l$.
Next, choose a nonnegative partition of unity 
$\,\phi_l; \, l=0,1,...\,$ 
subordinated to the cover 
$\,\mathcal D_l; \, l=0,1,...\,$ of $\mathcal D$
and also  choose, for each $l \geq 0$, nonnegative, smooth 
bump functions $\psi_l$, $0 \leq \psi_l \leq 1$, supported in 
$\mathcal D_l$ with $\psi_l \equiv 1$ on the support of $\phi_l$.
Then 
$\sum_{l \geq 0} \, \phi_l \, = 1$ and $ \psi_l \, \phi_l = \phi_l$ for all
 $l$.
\medskip

We aim to show that the operator
$M:\, C_{s,0}^{k+\alpha,p}(Q_\delta)
\to C_{s,0}^{k+2+\alpha,p}( Q_\delta)$
defined as
$$ M g \,= \, \sum_{} \psi_l M_l\, \phi_l g$$
satisfies 
$$||\, LM \,  g \,- \,g\, ||_{C_s^{k+\alpha,p} (Q_\delta)} < \frac 12 \,\,
||g||_{C_s^{k+\alpha,p} (Q_\delta)}  \qquad \forall g \in C_{s,0}^{k,\alpha,p}
(Q_\delta)$$
if the  cover $\{ \mathcal D_l \}$ and  $\delta$ are chosen appropriately.
Indeed, we can write 
$$ L\,M g - g  = \sum_{l} L \, \psi_l M_l\, \phi_l g -  \sum_{l} \phi_l g= 
\sum_{l} \,\psi_l \, (L M_l - I)\, \phi_l g  + 
\sum_{l}\,  [\,L,\psi_l \,]\, M_l\, \phi_l g$$
with $[\,L,\psi_l \,]$ denoting the commutator of $L$ and $\psi_l$.
The commutator  $[\, L, \psi_l \,]$ is only of first order
and it can be estimated as 
$$ || [\,L,\psi_l \,] \, M_l\phi_l g ||_{C_s^{k+\alpha,p} (Q_\delta)} \leq 
C \, \left ( \,
 || \vartheta\, D (M_l\phi_l g) ||_{C_s^{k+\alpha,p} (Q_\delta)} + 
\, ||M_l\phi_l g||_{C_s^{k+\alpha,p} (Q_\delta)}\, \right  ).$$
Let $\epsilon >0$.
%Since the function $\vartheta$ is proportional to the distance $d$
%to the boundary of $\mathcal D$, 
It follows via the  H\"older spaces interpolation from Lemma~\ref{holder} that 
$$|| \vartheta\, D (M_l\phi_l g) ||_{C^{k+\alpha,p}_s(Q_\delta)}
\leq \epsilon \,  ||M_l\phi_l g||_{C^{k+2+\alpha,p}_s (Q_\delta)} + C(\epsilon) 
\,  ||M_l\phi_l g||_{C^{k,p}(Q_\delta)}.$$
However, for each $k$ we have 
$$ ||M_l\phi_l\ g||_{C^{k+2+\alpha,p}_s (Q_\delta)} \leq C \, ||g||_{C_s^{k,
\alpha,p} (Q_\delta)}$$
and therefore, since $M_l\phi_l g \equiv 0$ at $t=0$,
$$||M_l\phi_l g||_{C^{k,p}_s(Q_\delta)} 
\leq C \, \delta\, ||g||_{C^{k+\alpha,p}_s (Q_\delta)}.$$
It follows that if we choose  $\delta$  sufficiently small:
$$\sum_{l} \, ||\, [\,L,\psi_l \,] \, M_l\phi_l g\, ||_{C^{k+\alpha,p}_s
 (Q_\delta)} \leq \frac 14 \,\,  ||g||_{C^{k+\alpha,p}_s(Q_\delta)}.$$
On the other hand we have $(\,L M_0 - I \,)\varphi_0 g=0$,
 while for $l \geq 1$,  we can make the norm of each
of the operators $L M_l - I$ arbitrarily  close to zero by
choosing the diameters of the domains  $\mathcal D_l$  sufficiently 
small:

$$ ||\, \sum_{l} \psi_l \, (L\,M_l - I)\, \phi_l g \, ||_{C^{k+\alpha,p}_s
(Q_\delta)} < \frac 14\,\,  ||g||_{C^{k+\alpha,p}_s(Q_\delta)}$$
for all $g \in C_s^{k+\alpha,p} (Q_\delta)$, if $\rho$ and $\delta$ are
both sufficiently small.
The above inequalities give:
$$ ||\, L\,M g - g\, || _{C^{k+\alpha,p}_s(Q_\delta)}
\leq \frac 12  \,\, ||g||_{C^{k+\alpha,p}_s(Q_\delta)}$$
for all $g \in 
C^{k+\alpha,p}_{s,0}(Q_\delta)$.
We conclude that for every $g \in C^{k+\alpha,p}_{s,0}(Q_\delta)$ there exists
a function $w \in C^{k+2+\alpha}_{s,0}(Q_\delta)$ such that  $Lw=g$.
In addition 
\begin{equation}\label{ineq:w}
 ||w||_{C^{k+2+\alpha,p}_s(Q_\delta)} 
\, \leq \, C \, ||g||_{C^{k+\alpha,p}_s(Q_\delta)}
\end{equation}
with $C$ depending only on $\mathcal D$ 
and the constants $\alpha$, $k$, $\lambda$ and $T$.

The last inequality implies we can extend the solution on a bigger interval. Hence, one can show that 
$$||w||_{C^{k+2+\alpha,p}_s(Q)} \leq C(T) \left ( ||w_0||_{C^{k+2+\alpha,p}_s(\mathcal{D})} + ||g||_{C^{k+\alpha,p}_s(Q)} \right )$$
where the constant $C(T)$ depends  only on the 
numbers $\alpha$, $k$, $\lambda$ and $T$.

\no This last inequality implies uniqueness.
\end{proof}

Finally, the following existence result follows from Theorem \ref{thm:exi2} and the Inverse
Function Theorem between Banach spaces along the line  of the proof of Theorem 7.3 in
\cite{DH2}:

\begin{theorem}\label{thm:exi3}
Let $w_0$ be a function in   $C^{k+2+\alpha,p}_s(\mathcal D)$.
 Assume that the linearization $DM(\bar w)$ of the 
fully-nonlinear operator
$$Mw = w_t - F(t,u,v,w,Dw, D^2w)$$
defined on the cylinder $Q = \mathcal D  \times [0,T]$,
satisfies the hypotheses of Theorem~\ref{thm:exi2}  at
all points $\bar w \in C^{k+2+\alpha,p}_s(Q)$, with 
$||\bar w - w_0||_{ C^{k+2+\alpha,p}_s(Q)} \leq \mu$, for some
$\mu >0$.
Then,  there exists a number $\tau_0$ in $0<\tau_0 \leq T$
depending on the  constants $\alpha,\,p$, $k$, $\lambda$
and $\mu$, for which the initial value problem
$$\left\{\begin{array}{ll} w_t = F(t,u,v,w,Dw, D^2w) 
\qquad &\text{in $\,\,\mathcal D \times [0,\tau_0]$}\\
w(\cdot,0)=w_0 \qquad &\text{on $\,\, \mathcal D$} \end{array}\right.$$
 admits a solution $w$ in the space $C^{k+2+\alpha,p}_s
(\mathcal D \times [0,\tau_0])$. Moreover,
$$||w||_{C^{k+2+\alpha,p}_s(\mathcal D \times [0,\tau_0])}
\leq C\, ||w_0||_{C^{k+2+\alpha,p}_s(\mathcal D)}$$
for some positive constant $C$ which depends only on $\alpha$, $p$, $k$,
$\lambda$ and $\mu$.
\end{theorem}

\section{Global change of coordinates and existence in $C^{k+2+\alpha,p}_s$}\label{global}

In this section we introduce a global change of coordinates which transforms the HMCF for a surface $\Sigma_0$ into a fully-nonlinear degenerate parabolic PDE on $\mathcal D$.

 Let $S$ be a smooth surface close to $\Sigma_2$. Let $S : \mathcal D \to \mathbb{R}^3$, indicate a
parameterization of $S$ on the unit disk $\mathcal D$.
 $$S(u,v)=(x,y,z)\in \R^3\,\,\,\hbox{ which maps }\,\,\,\partial \mathcal D\,\,\, \hbox{ onto}\,\,\, S\cap \{x=0\}$$ 
%The coordinates $x$ and $y$ are chosen to be smooth on $ \mathcal D$. 
\noindent We denote by $x_u, x_v, x_{uu}, x_{uv}, x_{vv}$ the partial derivatives of $x$ with respect to $u$ and $v$. The same notation will be used for the function $y$.
%Moreover, we require that:
%$$x_u, \,y_u\sim \vartheta^{1-p},\,
%|x_{uv}|\leq C\,\vartheta^{1-p}\,\,\,\,\hbox{and}\,\,\,\, x_{uu}\sim \vartheta^{2-p}$$
%\noindent where $\vartheta $ behaves like distance to the boundary as in Section~\ref{deg}. 

%\no and  we denote $x_u,x_v,x_w,y_u, y_v,y_w, z_u, z_v, z_w$ to be the partial
%derivatives of the functions $x =x(u,v,w), y = y(u,v,w)$ with respect to
%$u,v,w$. We use similar notation for the second derivatives of these
%functions.
\noindent Let $\eta >0$ be sufficiently small. Let $T = (T_1,T_2,T_3)$ be a
smooth vector field transverse to $S$.  We define the global change of
coordinates $\Phi : \mathcal D \times [-\eta,\eta] \to \mathbb{R}^3$ by

\begin{align}
\begin{pmatrix}
x\\
y\\
z	
\end{pmatrix}
= \Phi
\begin{pmatrix}
u\\
v\\
w
\end{pmatrix}
= S(u,v) + w T(u,v)
\end{align}

\no or more explicitly

\begin{align}
x &= S_1(u,v) + w T_1(u,v)\\
y &= S_2(u,v) + w T_2(u,v)\\
z &=S_3(u,v)+w T_3(u,v)
\end{align}

\no Let $\delta >0$ be small number, such that
$$T_3 (u,v) =0 \qquad \text{on }\,\, {\mathcal D} \setminus {\mathcal D}_{1-\delta}$$

\noindent denoting, as above,  the transverse vector field to the surface $\mathcal S$.
Notice that by choosing the smooth surface sufficiently close to 
the surface $z=h(x,y)$, we can make $\delta$ to depend only on 
the constant $\lambda$ which depends on the initial non-degeneracy conditions on the surface $\Sigma_0$.

%To show that $w$  satisfies  an equation of the desired form in
%the interior cylinder ${\mathcal D}_{1-\delta}\times [0,T]$ is straight forward.
%Hence, we will analyze ${\mathcal D}\setminus 
%{\mathcal D}_{1-\delta}$.
We write the 
 first and second derivatives
of $z$ with respect to $x,y$ and $t$  in terms of 
the first and second derivatives
of $w$ with respect to $u,v$ and $t$.

\noindent 
If $z= h_0(x,y)$ then we compute the first and second partial derivatives of $z$ with respect to $x$ and $y$ in terms of $w = l(u,v)$ seen as functions of $u$ and $v$.

%% AMS-LaTeX Created by Wolfram Mathematica 5.1

Let $A$ be the Jacobian matrix relative to the transformation of coordinates:
$$A  = \left (
\begin{split}
\frac{\partial u}{\partial x}\,\,\, & \frac{\partial v}{\partial x}\\
\frac{\partial u}{\partial y}\,\,\, & \frac{\partial v}{\partial y}
\end{split}\right )=
\left (
\begin{split}
a\,\,\, & b\\
c\,\,\, & d
\end{split}\right )
$$
Let $\nabla z$, $\nabla u$, $\nabla v$  be, respectively  the gradients of  $z$, $u$ and $v$:
$$\nabla z=\left(\begin{smallmatrix} \frac{\partial z}{\partial x}\\ \frac{\partial z}{\partial y}\end{smallmatrix}\right)=\left(\begin{smallmatrix} z_1\\ z_2\end{smallmatrix}\right)\,\,\,\nabla u=\left(\begin{smallmatrix} \frac{\partial u}{\partial x}\\ \frac{\partial u}{\partial y}\end{smallmatrix}\right)=\left(\begin{smallmatrix} u_1\\ u_2\end{smallmatrix}\right)\,\,
\nabla v=\left(\begin{smallmatrix} \frac{\partial v}{\partial x}\\
\frac{\partial v}{\partial
y}\end{smallmatrix}\right)=\left(\begin{smallmatrix} v_1\\
v_2\end{smallmatrix}\right)$$
\noindent We denote by $D^2 u$ and $D^2 v$ the following matrices:

$$D^2 u  = \left (
\begin{split}
\frac{\partial^2 u}{\partial x^2}\,\,\, & \frac{\partial^2 u}{\partial x\partial y}\\
\frac{\partial^2 u}{\partial x\partial y}\,\,\, & \frac{\partial^2 u}{\partial y^2}
\end{split}\right )=
\left (
\begin{split}
a_1\,\,\, & c_1\\
c_1\,\,\, & c_2
\end{split}\right )
$$

$$D^2 v  = \left (
\begin{split}
\frac{\partial^2 v}{\partial x^2}\,\,\,& \frac{\partial^2 v}{\partial x\partial y}\\
\frac{\partial^2 v}{\partial x\partial y}\,\,\,  &\frac{\partial^2 v}{\partial y^2}
\end{split}\right )=
\left (
\begin{split}
b_1\,\,\, & d_1\\
d_1\,\,\, & d_2
\end{split}\right )
$$

\noindent Based on the above, we define $a_2:=c_1\,$, $b_2:=d_1$ and denote 
$e_1=(1,0)\,\,\,e_2=(0,1)$ to be the basis vectors. 
Let $A^{-1}$ be the inverse matrix of $A$:
$$A^{-1}: = \left (
\begin{split}
x_u,\,\, & x_v\\
y_u\,\,\, & y_v
\end{split}\right )=
\left (
\begin{split}
x_1\,\,\, & x_2\\
y_1\,\,\, & y_2
\end{split}\right )
$$

\noindent Next, we introduce the matrices $B_1$ and $B_2$ which denote,
respectively, the derivative of the inverse matrix of $A$ and $A^{-1}$
with respect to $x$ and $y$;  $B_1=\frac{\partial A^{-1}}{\partial x}$,
$B_2=\frac{\partial A^{-1}}{\partial y}$ which can be computed as:

$$B_1 = \left (
\begin{split}
a\,x_{11}+b\,x_{12}\,\,\,\,\,\,& a\,y_{11}+b\,y_{12}\\ 
a\,x_{12}+b\,x_{22}\,\,\,\,\,\,& a\,y_{12}+b\,y_{22}
\end{split}\right )$$

$$B_2 = \left (
\begin{split}
c\,x_{11}+d\,x_{12}\,\,\,\,\,\,& c\,y_{11}+d\,y_{12}\\ 
c\,x_{12}+d\,x_{22}\,\,\,\,\,\,& c\,y_{12}+d\,y_{22}
\end{split}\right )$$

\noindent Note that we are using the following notation:
$$x_{11}=x_{uu},\, x_{12}=x_{uv},\,x_{22}=x_{vv},\, y_{11}=y_{uu},\,
y_{12}=y_{uv},\,y_{22}=y_{vv}.$$ \noindent The coefficients of the matrices $D^2 u$ and $D^2 v$ are evaluated as follows:
$$
\begin{array}{lcl}
\left(\begin{smallmatrix} a_1\\
c_1\end{smallmatrix}\right) &  =  -A\cdot B_1\cdot \nabla u;\,\,
\left(\begin{smallmatrix} b_1\\
d_1\end{smallmatrix}\right) & = -A\cdot B_1\cdot \nabla v\\
\\
\left(\begin{smallmatrix} a_2\\
c_2\end{smallmatrix}\right) & =-A\cdot B_2\cdot \nabla u;\,\,
\left(\begin{smallmatrix} b_2\\
c_2\end{smallmatrix}\right) & =-A\cdot B_2\cdot \nabla v
\end{array}
$$

\noindent Since: $\qquad \nabla z=A\cdot \left(\begin{smallmatrix} z_u\\
z_v\end{smallmatrix}\right)+z_w\left(\begin{smallmatrix}
\,\frac{\partial w}{\partial x}\\ \frac{\partial w}{\partial
y}\end{smallmatrix}\right)$, we obtain that:

$$\frac{\partial}{\partial x}\nabla z=A_x\cdot
\left(\begin{smallmatrix} z_u\\ z_v\end{smallmatrix}\right)+A\cdot
\left(\begin{smallmatrix} a\,z_{uu}+b\,z_{uv}+z_{uw}\frac{\partial
w}{\partial x}\\a\,z_{uv}+b\,z_{vv}+z_{vw}\frac{\partial w}{\partial
x}\end{smallmatrix}\right) +
(a\,z_{uw}+b\,z_{vw})\left(\begin{smallmatrix} \,\frac{\partial
w}{\partial x}\\ \frac{\partial w}{\partial
y}\end{smallmatrix}\right)+z_w\,\left(\begin{smallmatrix}
\,\frac{\partial^2 w}{\partial^2 x}\\ \frac{\partial^2 w}{\partial
x\partial y}\end{smallmatrix}\right).$$

$$\frac{\partial}{\partial y}\nabla z=A_y\cdot
\left(\begin{smallmatrix} z_u\\ z_v\end{smallmatrix}\right)+A\cdot
\left(\begin{smallmatrix} c\,z_{uu}+d\,z_{uv}+z_{uw}\frac{\partial
w}{\partial y}\\ c\,z_{uv}+d\,z_{vv}+z_{vw}\frac{\partial w}{\partial
y}\end{smallmatrix}\right) +
(c\,z_{uw}+d\,z_{vw})\left(\begin{smallmatrix} \,\frac{\partial
w}{\partial x}\\ \frac{\partial w} {\partial
y}\end{smallmatrix}\right)+z_w\,\left(\begin{smallmatrix}
\,\frac{\partial^2 w}{\partial x\partial y}\\ \frac{\partial^2 w}
{\partial^2 y}\end{smallmatrix}\right). $$

\noindent The gradient of the function $w$, as well as its partial
derivatives, can be expressed by using the matrix $A$:

$$
\left(\begin{smallmatrix} \,\frac{\partial w}{\partial x}\\ \frac{\partial w}{\partial y}\end{smallmatrix}\right)=
\left(\begin{smallmatrix} \,a\,w_u+b\,w_v\\ c\,w_u+d\,w_v\end{smallmatrix}\right)$$
$$
\left(\begin{smallmatrix} \,\frac{\partial^2 w}{\partial^2 x}\\ \frac{\partial^2 w}{\partial x\partial y}\end{smallmatrix}\right)=
\left(\begin{smallmatrix} \,a_1\,w_u+a\,(a\,w_{uu}+b\,w_{uv})+b_1\,w_v+b(a\,w_{uv}+b\,w_{vv})
\\ c_1\,w_u+c(a\,w_{uu}+b\,w_{uv})+d_1\,w_v+d(a\,w_{uv}+b\,w_{vv})\end{smallmatrix}\right)$$

$$
\left(\begin{smallmatrix} \,\frac{\partial^2 w}{\partial x\partial y}\\ \frac{\partial^2 w}{\partial^2 y}\end{smallmatrix}\right)=
\left(\begin{smallmatrix} \,a_2\,w_u+a\,(c\,w_{uu}+d\,w_{uv})+b_2\,w_v+b(c\,w_{uv}+d\,w_{vv})
\\ c_2\,w_u+c(c\,w_{uu}+d\,w_{uv})+d_2\,w_v+d(c\,w_{uv}+d\,w_{vv})\end{smallmatrix}\right).$$

\noindent 

Using the substitutions above we get:
$$\frac{\partial^2 z}{\partial^2
x}=A^1_{11}\,w_{11}+A^1_{12}\,w_{12}+A^1_{22}\,w_{22}+B^1_1\,w_1+B^1_2\,w_2+B^1_{12}\,w_1\,w_2+
B_{11}^1\,w_1^2+B_{22}^1\,w_2^2+
C_1$$
$$\frac{\partial^2 z}{\partial^2
y}=A^2_{11}\,w_{11}+A^2_{12}\,w_{12}+A^2_{22}\,w_{22}+B^2_1\,w_1+B^2_2\,w_2+B^2_{12}\,w_1\,w_2+
B_{11}^2\,w_1^2+B_{22}^2\,w_2^2+
C_2$$
$$\frac{\partial^2 z}{\partial x\partial y
}=A^{\circ}_{11}\,w_{11}+A^{\circ}_{12}\,w_{12}+A^{\circ}_{22}\,w_{22}+B^{\circ}_1\,w_1+B^{\circ}_2\,w_2+B^{\circ}_{12}\,w_1\,w_2+
B_{11}^{\circ}\,w_1^2+B_{22}^{\circ}\,w_2^2+
C_{\circ}$$

\noindent where the coefficients $A^{k}_{i,j},\,B_i^j,\,C^i$ with $k=0,1,2$, $i,j=1,2$ are defined as follows:
$$
\begin{array}{lcl}
A_{11}^1: & = & -a^2\,(\,-z_{w}+b\,y_w(z_u+w_1\,z_w)+d\,y_w(z_v+w_2 z_w))\\

A_{22}^1: & = & -b^2\,(\,-z_{w}+b\,y_w(z_u+w_1\,z_w)+d\,y_w(z_v+w_2 z_w))\\

A_{12}^1: & = & -2\,a\,b\,(\,-z_{w}+b\,y_w(z_u+w_1\,z_w)+d\,y_w(z_v+w_2 z_w))\\

A_{11}^2: & = & c^2(z_w - d\,y_w (z_v + w_2\,z_w))\\
A_{22}^2: & = & d^2(z_w - d\,y_w (z_v + w_2\, z_w))\\

A_{12}^2: & = & -2 c d (-z_w + d\, y_w (z_v + w_2 z_w))
 \end{array}
 $$

$$
\begin{array}{lcl}
A_{11}^0: & = & -a c (-z_w + b y_w (z_u + w_1 z_w) + d y_w (z_v + w_2 z_w))\\

A_{22}^0: & = & -b d (-z_w + b y_w (z_u + w_1 z_w) + d y_w (z_v + w_2 z_w))\\

A_{12}^0: & = & -b (2 a c y_{uw} + b c y_{vw} + a d y_{vw}) z_w\\
B_{11}^0: & = & -b (2 a c y_{uw} + b c y_{vw} + a d y_{vw}) z_w\\

B_{22}^0: & = & -d (b c y_{uw} + a d y_{uw} + 2 b d y_{vw}) z_w\\

B_{12}^0 : & = &-((b^2 c + a (b + 2 c)  d) y_{uw} + d(b(2 b + c) + a d) y_{vw}) z_w \\

B_1^0 : & = & -a^2(c x_{uu} + d x_{uv}) z_w - b(b c y_{vw} z_u + c d y_{vw} z_v -c z_{vw} + b c y_{uv} z_w +\\
{} & {} & b d y_{vv} z_w) - a (-2 c z_{uw} + 2 c d y_{uw} z_v + d^2 y_{vw} z_v - d z_{vw} + b (2 c y_{uw} z_u + d y_{vw} z_u +\\
{} & {} &  c (x_{uv} + y_{uu}) z_w + d (x_{uv} + y_{uv}) z_w))\\

B_2^0   : & = & -b^2 (c y_{uw} + 2 d y_{vw})z_u - b (-c z_{uw} + d (a y_{uw} z_u + c y_{uw} z_v + 2 d y_{vw} z_v - 2 z_{vw}) +\\
{} & {} & (c(c + d) x_{uv} + c d y_{uv} + d^2 y_{vv}) z_w) - a (c^2 x_{uu} z_w + d (-z_{uw} + c (x_{uv} + y_{uu}) z_w +d(y_{uw} z_v + y_{uv} z_w)))\\

C^0 : & = & -a^2 (c x_{uu} + d x_{uv}) z_u - a (b (c (x_{uv} + y_{u}u) + d ( x_{uv} + y_{uv})) z_u - c z_{uu} + c^2 x_{uu} z_v + \\
{} & {} &c d (x_{uv} + y_{uu})  z_v + d (-z_{uv} + d y_{uv} z_v)) - b (b (c  y_{uv} + d y_{vv}) z_u - c z_{uv} + c^2 x_{uv}z_v + c \\
 {} & {} &         d (x_{uv} + y_{uv})z_v + d (d y_{vv}z_v - z_{vv}))

 \end{array}
$$
$$
\begin{array}{lcl}
B_{11}^1: & = & -2 a b ( a y_{uw} + b y_{vw})\\

B_{22}^1: & = & -2 b d ( a y_{uw} + b y_{vw})\\

B_{12}^1   : & = &  -2 ( b^2 + a d ) ( a y_{uw} + b y_{vw})\,z_w\\
B_{11}^2: & = & 0\\
B_{22}^2: & = & -2\, d^2(c\, y_{uw} + d\, y_{vw})\, z_w\\
B_{12}^2    : & = &  -2 c d (c y_{uw} + d y_{vw}) z_w\\

B_1^1  : & = & -2 a(a(b\,y_{uw}\, z_u - z_{uw} + d\,y_{uw} z_v) + b(b\, y_{vw}\, z_u + d\, y_{vw}\, z_v - z_{vw})) - \\
{} & {} & (a^3\, x_{uu}  +  a^2 b(2\, x_{uv} + y_{uu}) + a b^2(x_{uv} + 2\, y_{uv}) + b^3 y_{vv}) z_w\\
B_2^1 : & = & -a^2(c x_{uu} + d\ y_{uu})z_w - 2 a b(b y_{uw}z_u - z_{uw} + d y_{uw}z_v + c x_{uv} z_w + d y_{uv} z_w) - \\
{} & {} & b^2  (  2 b y_{vw} z_u + 2 d y_{vw} z_v - 2 z_{vw} + c x_{uv} z_w + d y_{vv} z_w)
 \end{array}
$$

$$
\begin{array}{lcl}
C_1: & = & -a^3 x_{uu} z_u - a^2(b(2 x_{uv} + y_{uu}) z_u - z_{uu} + c x_{uu} z_v + d y_{uu}z_v) - \\
{} & {} & a b(b (x_{uv} + 2 y_{uv}) z_u - 2 z_{uv} + 2(c x_{uv} + d y_{uv}) z_v) - b^2(b y_{vv}z_u + c x_{uv} z_v + d y_{vv} z_v - z_{vv})\\

B_1^2: & = & 2 c((c(z_{uw} - d y_{uw} z_v) + d(-d y_{vw} z_v + z_{vw})) -  a d^2 x_{uv} z_w\\
B_2^2 :& = &-2 d(-c z_{uw} + c d y_{uw} z_v + d^2 y_{vw} z_v - d z_{vw}) - (c^3 x_{uu} + c^2\\
 {} & {} &  d(2 x_{uv} + y_{uu}) + c d^2(x_{uv} + 2 y_{uv}) + d^3y_{vv}) z_w\\
 C_2 : & = & -a d^2 x_{uv} z_u - c^3 x_{uu} z_v + c^2(z_{uu} - d (2 x_{uv} + y_{uu}) z_v)+ \\
 {} & {} &  c d(2 z_{uv} - d (x_{uv} + 2 y_{uv})z_v) + d^2(-d y_{vv} z_v + z_{vv}) 

\end{array}
$$

\newpage

\no {\it Evolution of $w_t$}. 
\no
The evolution equation of $w$ is :
\smallskip
\begin{equation}\label{eqnwt}
w_t=\displaystyle{\frac{1}{z_y\,y_w-z_w}\,z_t}
\end{equation}
\smallskip
\no where the function $z$ satisfies the non linear PDE:
\begin{equation} 
z_t=
\frac{z_{xx}z_{yy}-z_{xy}^2 }{(1+z_y^2)z_{xx}-2\,z_x z_y
z_{xy}-(1+z_x^2)z_{yy}}
\end{equation}

%\noindent On the other hand, after some easy but long computation we obtain:
% $$\frac{\partial^2 z}{\partial^2 x}\,\frac{\partial^2 z}{\partial^2
% y}\,-\left(\frac{\partial^2 z}{\partial x\partial y}\right
%)^2\,=$$

%$$\gamma_{11}\,w_{11}^2+2\,\gamma_{12}\,w_{12}^2+\gamma_{22}\,w_{22}^2+\rho_{11}\,w_{11}+\rho_{12}\,w_{12}+\rho_{22}\,w_{22}+\gamma_{1}\,w_1^2+\gamma_{2}\,w_2^2+\rho_{1}\,w_1+\rho_2\,w_2+\Theta
% $$
%%\noindent while the denominator:

\noindent By replacing the partial derivatives of $z$ in terms of the derivatives of $w$, we find the linearization $L$ to the fully non linear operator given by the Equation~\ref{eqnwt}.

\no {\it  Linearization}. Let $\tilde w$ be a point close to the initial data $w_0$. Then, 
%We assume $w_2,\, z_v,\,z_{vw},\, x_w,\,z_{vw},\, x_{vv},\,x_w,\,x_{vw},\,x_{uw}\sim 0$. 
the linearization of the operator given by the Equation~\ref{eqnwt} is given by:

$$L(\tilde{w})=a_{11}\,\tilde{w}_{11}+2\,a_{12}\,\tilde{w}_{12}+a_{22}\,\tilde{w}_{22}+b_1\,\tilde{w}_1+b_2\,\tilde{w}_2+c\,\tilde{w}+d $$

\noindent where its coefficients behave as follows:
$$  
\begin{array}{ll}
& a_{11}\approx \frac{ w_1^2}{w_{11}^2}\, g_{11};\,\, a_{12}\approx \frac{ w_{11}\,w_1}{w_{11}^2}\,g_{12}\\
& a_{11}\approx \frac{ w_{11}^2}{w_{11}^2}\,g_{22};\,\, b_{1}\approx \frac{ w_{11}\,w_1}{w_{11}^2}\,h_1\\
& b_{2}\approx \frac{ w_{11}^2}{w_{11}^2}\,h_2
\end{array}
$$
\noindent where $\{g_{ij}\}$, $h_i$, $c$ and $d$ with $i,j=1,2$ are functions which belong to the space $C^{\alpha}_s$. In particular $g_{22}$ and $ b_2$ belong to $C^{\alpha,p}_s$. It follows that the coefficients of the operator $L$ satisfy the same conditions of the operator of Theorem~\ref{thm:exi2}. 
\smallskip

\noindent {\it We say that $\Sigma_0$ is of  H\"older class $C^{k+2+\alpha,p}_s$ if the function $w_0$ obtained belongs 
the space $C^{k+2+\alpha,p}_s(\mathcal D)$.}
\smallskip

 Based on the above definition, we state the following theorem:

\begin{theorem}\label{thm:exi4}
Assume that  the  initial surface $\Sigma_0$ belongs to the class
$C^{k+2+\alpha,p}_s$ and satisfies the non-degeneracy condition~($\star$) .  
Then, under the coordinate change $\Phi$,  the HMCF with initial data the surface $\Sigma_0$
converts into  the  initial value problem
\[
\left\{
\begin{array}{ll}
 Mw =0
& on\,\,{\mathcal{D}}\times [0,T]\\
w=w_0 & at \,\,\,\,\,\,\,t=0%
\end{array}%
\right.
\]
with $w_0 \in  C^{k+2+\alpha,p}_s(\mathcal D)$
and 
$$Mw= w_t - F(t,u,v,w,Dw, D^2w)$$
satisfying the hypotheses of Theorem \ref{thm:exi3}.
\end{theorem}
\no As an immediate Corollary of Theorem \ref{thm:exi3} and Theorem \ref{thm:exi4}
we obtain the following existence result:

\begin{thm}
Under the same assumptions as in Theorem \ref{thm:exi4}, there 
exists a number $\tau_k >0$ for which the HMCF
with initial data the surface $\Sigma_0$ admits a solution $\Sigma_t$
on $0 \leq t \leq \tau_k$. Moreover, under the
coordinate change $\Phi$  the strictly convex part
$\Sigma_2(t)$ of $\Sigma_t$ is converted  to a function $w(t)$
which belongs to the H\"older class $C^{k+2+\alpha,p}_s(Q_k)$,
on $Q_k = \mathcal D \times [0,\tau_k]$.
\end{thm}

%In the next section we will actually show that the
%function $w$ is smooth up to the boundary of
%$\mathcal D$ for all $t \in  (0,T]$, for 
%some $T >0$. 
%\bigskip

\begin{theorem}\label{thmreg}
Assume that the initial surface $\Sigma_0$ satisfies the assumptions of
Theorem~\ref{thm1}. Then, the solution $\Sigma_t$ of the HMCF is converted, via
the coordinate change studied in Section~\ref{global}, to a function
$w$ which belongs, for any positive integer $k$, to the H\"older
class $C^{k+2+\alpha,p}_s(Q)$, on $Q = \mathcal D \times
(0,T]$. Moreover, for any $\tau$ in $0<\tau<T$ we have
\begin{equation}\label{eqn-w}
||w||_{ C^{k+2+\alpha,p}_s(\mathcal D \times [\tau,T])}
\leq C_k(\tau, ||w_0||_{C^{2+\alpha,p}_s(\mathcal D)})
\end{equation}
\end{theorem}
%\bigskip
\begin{proof}
We omit the proof as it is similar to the one done in~ \cite{DH2}. Also, for more details we invite the reader to look at the Ph.D. thesis~\cite{Cap} of the first author.
\end{proof}
\section{The  proof of the Main Theorem}\label{sec-main}

In this section we will give the proof of the Main Theorem
stated in the Introduction. We will actually prove the following stronger result, where
we relax the regularity assumptions on the initial surface.

\begin{theorem}\label{thm6}
 Assume that the strictly convex part $\Sigma_2$ of
the initial surface $\Sigma_0$ belongs to the class
$C^{2+\alpha,p}_s$ and satisfies the 
non-degeneracy conditions~($\star$). Then, the {\it HMCF}
$$\frac {\partial P}{\partial t} =\frac{ K}{H}\, \stackrel{\rightarrow}{N}
         \quad \qquad  t  \in  [0,T]$$
with initial data the surface $\Sigma_0$
admits a solution $\Sigma_t$ which is smooth up to the interface,
for $0<t\leq T$. In particular, the interface $\Gamma_t$ is a 
smooth curve for every $0<t\leq T$ which moves by the curve shortening Flow.
\end{theorem}

\noindent{\em Remark.} It can be easily checked that if the initial surface satisfies the
conditions of the Main Theorem, then it will satisfy the weaker
conditions of Theorem~\ref{thm6}.

\begin{proof}
 Assume that the strictly convex part
$\Sigma_2$ of the initial surface $\Sigma_0$ belongs to the class
$C^{2+\alpha,p}_s$ and satisfies the non-degeneracy condition~($\star$), then we have proven existence for the
{\it HMCF} in Theorem~\ref{thm:exi4}.\\
\no From Theorem~\ref{thm:exi4} we have that $w \in C^{k+2+\alpha,p}_s(\mathcal D
\times (0,T])$,  for all nonnegative integers $k$.
In particular this implies that for all 
integers $k$ we have $w(t) \in C^{k,p}(\mathcal D$, 
 for all $\tau$ in $0<t<T$.
It follows that  $w$ is $C^{\infty,p}$ smooth up to the boundary of $\mathcal D$.
Going back to the original coordinates, we conclude that the
strictly convex part of the surface $\Sigma_t$, $0 < t \leq T$ is smooth 
up to $z=0$  and that the interface $\Gamma_t$ is smooth.

\end{proof}

\section{Comparison principle}\label{sec-cp}
In this final section we will give the proof of the comparison principle
for the HMCF and we will show that the solution given in the Main Theorem is a viscosity solution.

\begin{prop} {\bf (Comparison principle)}
Let $\Sigma_0$ be a surface of class $C^{2+\alpha,p}_s$ that 
 satisfies condition~($\star$), and let $\Sigma^{+}$ be a smooth, strictly
 convex surface containing $\Sigma_0$ at time $t=0$, then the surface
 $\Sigma_t$ obtained by evolving $\Sigma_0$ by the HMCF,
 is contained in the surface $\Sigma^{+}_t$ obtained by evolving
 $\Sigma^{+}_0$ by the HMCF up to the time of existence
 of $\Sigma_t$.  Analogously, if $\Sigma_0$ contains a smooth,
 strictly convex surface $\Sigma^{-}$ at time $t=0$, then the surface
 $\Sigma_t$ contains the surface $\Sigma^{-}_t$ obtained by evolving
 $\Sigma^{-}_0$ by the HMCF up to the time of
 existence of $\Sigma_t$.

\end{prop}

\begin{proof} We observe that by the the classical  maximum principle the surfaces $\Sigma_t$ and $\Sigma^{+}_t$ cannot touch were they are both strictly convex. 
Next, we suppose that there exists a time $\bar{t}$ where they first touch at a
point $\bar{P}$, then this cannot happen in the interior of
the flat side. Hence $\bar{P}$ belongs to the boundary of the
flat side. Suppose $\Sigma_t$ has the flat side on the $x=0$ plane,
then the tangent plane to the surface at the point $\bar P$ would be
contained in the $x=0$ plane. Now if the two regions touch, then, because they are of
class $C^1$, the tangent to $\Sigma^{+}_{\bar{t}}$
at $\bar P$ would be contained in the $x=0$ plane. But
$\Sigma_{\bar{t}}^{+}$ is strictly convex, hence this would imply that
a part of $\Sigma_{\bar{t}}^{+}$ is inside $\Sigma_{\bar{t}}$ which
leads to a contradiction. The second part of the proof is straightforward since if
$\Sigma$ contains a smooth, strictly convex surface $\Sigma^{-}$ at
time $t=0$, then the two surfaces cannot touch at the flat side of
$\Sigma_t$ because the flat side does not move in its normal
direction. Once again, by the classical maximum principle for parabolic
equations the two surfaces cannot touch where they are strictly convex
either. 
\end{proof}

\begin{corollary}{\bf (Viscosity solutions)} Let $h\in \mathbb{N}$, $0<\gamma\leq 1$.
Let  $\Sigma_0$ be a convex surface of class $C^{h,\gamma}$ that satisfies the hypothesis of the Main Theorem. Then the solution $\Sigma_t$ given by Theorem~\ref{thm6} is a viscosity solution of class $C^{h,\gamma}$. Moreover, this solution is unique upon satisfying the non-degeneracy condition~$(\star)$.
\end{corollary}
\begin{proof} It is a consequence of the comparison principle. \end{proof}

%\section{Constant rate of decrease for the area}\labal{sec-area}
%Finally we observe the following property of surfaces evolving by the {\it harmonic mean curvature flow}. We hope that this property will have useful applications.
%  \begin{e-proposition}{\bf (Constant rate of decrease for the Area)}
%Let $A$ be the surface  area of a smooth, closed, strictly convex surface $\Sigma_0$
%evolving by HMCF. Then, we have

%\begin{equation}\label{eqn:area}
%\frac{dA}{dt}=-4\pi
%\end{equation}
%\end{e-proposition}
%\begin{proof} For a strictly convex surface evolving with a certain speed $v$ in the direction of the inward normal vector $ \stackrel{\rightarrow}{N}$ by the equation $\frac{d\,P}{d\,t}=v\,
%\stackrel{\rightarrow}{N}$ the area decreases by the First Variation
%Formula. This means that:
%\begin{equation}\label{eqn:areaint}
%\frac{dA}{dt}=-\int_{\Sigma} v\,H dp
%\end{equation}
%\no where $H$ indicates the mean curvature. For more details about how to derive this general formula see~\cite{ICh}.
%Thus, for the {\it HMCF}, $v=K/H$, which implies by the Gauss-Bonnet formula that the area shrinks at a constant rate:
%$$\frac{dA}{dt}=-\int_{\Sigma} v\,H dp=-\int_{\Sigma} K dp =-4\pi$$
%\end{proof}

\end{document}